\documentclass{article}
 \usepackage{amsmath}

\newtheorem{thm}{Theorem}[section]

\newtheorem{prop}[thm]{Proposition}

\newtheorem{lem}[thm]{Lemma}

\newtheorem{rem}[thm]{Remark}

\newtheorem{ex}[thm]{Example}

\newcommand{\be}{\begin{equation}}
\newcommand{\ee}{\end{equation}}
\newcommand{\ben}{\begin{enumerate}}
\newcommand{\een}{\end{enumerate}}
\newcommand{\beq}{\begin{eqnarray}}
\newcommand{\eeq}{\end{eqnarray}}
\newcommand{\beqn}{\begin{eqnarray*}}
\newcommand{\eeqn}{\end{eqnarray*}}

\newcommand{\pa}{\partial}

\newcommand{\pxi}{ {\pa \over \pa x^i}}

\newcommand{\pyi}{{\pa \over \pa y^i}}

\begin{document}
\title{A Note on a Class of   Finsler Metrics of Isotropic S-Curvature}
\author{Guojun Yang }
\date{}
\maketitle
\begin{abstract}
An $(\alpha,\beta)$-metric is defined by a Riemannian metric and
$1$-form.  In this paper, we investigate  the known
characterization for $(\alpha,\beta)$-metrics of isotropic
S-curvature. We show that such a characterization should hold in
dimension $n\ge 3$, and for the 2-dimensional case, there is one
more class of isotropic S-curvature than the higher dimensional
ones. Further, we construct corresponding examples for every
two-dimensional class, especially for the class that the norm of
$\beta$ with respect to $\alpha$ is not a constant.

\


\end{abstract}

\section{Introduction}
The S-curvature is one of the most important non-Riemannian
quantities in Finsler geometry which was originally introduced for
the volume comparison theorem (\cite{shen1}). Recent studies show
that the S-curvature plays a very important role in Finsler
geometry (cf. \cite{BRS}, \cite{CMS},
 \cite{shen4}--\cite{Y1}). It is proved
that, if an $n$-dimensional Finsler metric $F$ is of {\it
isotropic}
 S-curvature ${\bf S}=(n+1)c(x)F$ for a scalar function $c(x)$ and of
 scalar flag curvature ${\bf K}={\bf K}(x,y)$, then {\bf K} can be
 given by
 $$
 {\bf K}=\frac{3c_{x^m}y^m}{F}+\tau(x),
 $$
where $\tau(x)$ is a scalar function (\cite{CMS}).

An $(\alpha,\beta)$-metrics is a Finsler metrics defined by a
Riemann metric $\alpha=\sqrt{a_{ij}(x)y^iy^j}$ and a 1-form
$\beta=b_i(x)y^i$,
 $$F=\alpha \phi(s),\ \ s=\beta/\alpha,$$
 where $\phi(s)$ satisfies certain conditions such that $F$ is a regular Finsler
 metric. A special class of $(\alpha,\beta)$-metrics are Randers metrics defined by $F
 =\alpha+\beta$. By aid of navigation method, there have been
 given
 some important characterizations or local structures for  Randers metrics  of
 isotropic S-curvature (\cite{ShX}, \cite{Y}, \cite{Y1}).
 Especially, it is known that the local structure has been solved for two-dimensional Randers metrics  of
 isotropic S-curvature (\cite{Y}).

 In this paper, we show there is a gap in the known
characterization for $n(\ge 2)$-dimensional
$(\alpha,\beta)$-metrics of isotropic S-curvature (\cite{CS1}),
and there is one more new class  in $n=2$.

 For a pair of $\alpha$ and $\beta$, we define
 $b:=||\beta||_{\alpha}$ and
 $$r_{ij}:=\frac{1}{2}(b_{i|j}+b_{j|i}),\ \
 s_{ij}:=\frac{1}{2}(b_{i|j}-b_{j|i}),\ \
 r_j:=b^ir_{ij},\ \ s_j:=b^is_{ij},\ \ s^i:=a^{ik}s_k,$$
 where $b_{i|j}$ denote the covariant derivatives of $\beta$ with respect to
 $\alpha$ and
 $b^i:=a^{ij}b_j$, and $(a^{ij})$ is the inverse of
 $(a_{ij})$. For a $C^{\infty}$ function $\phi(s)>0$ on
 $(-b_o,b_o)$, define
  \be\label{y01}
 \Phi:=-(Q-sQ')(n\Delta+sQ+1))-(b^2-s^2)(1+sQ)Q'',
 \ee
 where $\Delta:=1+sQ+(b^2-s^2)Q'$ and $Q:=\phi'/(\phi-s\phi')$. It
 is known in \cite{CS} that a Randers metric $F=\alpha+\beta$ is of
 isotropic S-curvature ${\bf S}=(n+1)c(x)F$ if and only if
 $$r_{ij}=2c(a_{ij}-b_ib_j)-b_is_j-b_js_i.$$
In this paper we obtain the following theorem based on the result
in \cite{CS1}.

\begin{thm}\label{th1}
  Let $F=\alpha \phi(s)$, $s=\beta/\alpha$, be an $n$-dimensional
  $(\alpha,\beta)$-metric, where $\phi(0)=1$. Suppose that
   $\phi(s)\ne k_1\sqrt{1+k_2s^2}+k_3s$ for any constants $k_1>0, k_2$ and $k_3$. $F$ is  of isotropic S-curvature
  ${\bf S}=(n+1)c(x)F$ if and only if one of the following holds
 \ben
   \item[{\rm (i)}] {\rm($n\ge 2$)} $\beta$ satisfies
      \be\label{y1}
        r_j+s_j=0,
      \ee
    and $\phi=\phi(s)$ satisfies
     \be\label{y2}
      \Phi=0.
     \ee
     In this case, ${\bf S}=0$.
     \item[{\rm (ii)}] {\rm($n\ge 2$)} $\beta$ satisfies
     \be\label{y3}
       r_{ij}=\epsilon(b^2a_{ij}-b_ib_j),\ \ s_j=0,
      \ee
    and $\phi=\phi(s)$ satisfies
     \be\label{y4}
      \Phi=-2(n+1)k\frac{\phi\Delta^2}{b^2-s^2},
     \ee
     where $k$ is a constant. In this case, ${\bf S}=(n+1)cF$ with $c=k\epsilon$.
     \item[{\rm (iii)}] {\rm($n\ge 2$)} $\beta$ satisfies
     \be\label{y5}
       r_{ij}=0,\ \ s_j=0,
      \ee
    and $\phi=\phi(s)$ is arbitrary. In this case, ${\bf S}=0$.

    \item[{\rm (iv)}] {\rm($n= 2$)} $\beta$ satisfies
      \be\label{y6}
         r_{ij}=\frac{3k_1+k_2+4k_1k_2b^2}{4+(k_1+3k_2)b^2}(b_is_j+b_js_i),
      \ee
    and $\phi=\phi(s)$ is given by
     \be\label{y7}
      \phi(s)=\big\{(1+k_1s^2)(1+k_2s^2)\big\}^{\frac{1}{4}}e^{\int^s_0\tau(s)ds},
     \ee
     where $\tau(s)$ is defined by
      \be\label{y8}
       \tau(s):=\frac{\pm\sqrt{k_2-k_1}}
      {2(1+k_1s^2)\sqrt{1+k_2s^2}},
      \ee
      and $k_1$ and $k_2$ are constants with $k_2>k_1$. In this case, ${\bf S}=0.$
  \een

\end{thm}

We find a gap in the proof of Lemma 6.2 in the characterization
for $n$-dimensional $(\alpha,\beta)$-metrics of isotropic
S-curvature in \cite{CS1} and it should hold for $n\ge 3$ (see the
former three classes in Theorem \ref{th1}). In this paper,
starting from Lemma 6.1 in \cite{CS1}, we give a new proof and
then we obtain an additional class for $n=2$, that is Theorem
\ref{th1} (iv). The proof from (\ref{y60}) below to the end of
Section \ref{sec4} gives Theorem \ref{th1} (iv) in two-dimensional
case.

Note that in this paper, to complete the proof of Theorem
\ref{th1}, we only need to discuss the solutions of equations
(\ref{y17}), (\ref{y18}) and (\ref{y20}) below in the  case $b\ne
constant$ (see Section \ref{sec3} and Section \ref{sec4} below)
(the case $b=constant$ has been solved in \cite{CS1}). It turns
out that the case $b\ne constant$ occurs indeed if the dimension
$n=2$. Then combining with the discussions in \cite{CS1} we
finally give the complete proof of Theorem \ref{th1}.

Based on the class as shown in Theorem \ref{th1} (iv), we further
prove in \cite{Y3} that $F$ in Theorem \ref{th1} (iv) is
positively definite if and only if
 $$1+k_1b^2>0.$$
Then we  also prove in \cite{Y3}  that $b=constant$ in (\ref{y6})
if and only if
$$(1+k_1b^2)(1+k_2b^2)s_i=0,$$
and so by the positively definite condition $1+k_1b^2>0$, it is
easy to see that $s_i=0$ and thus $\alpha$ is flat and $\beta$ is
parallel with respect to $\alpha$.

In (\ref{y7}) and (\ref{y8}), it is easy to see that if $k_1=k_2$,
then $\phi(s)=\sqrt{1+k_1s^2}$. So the case $k_1=k_2$ is excluded.
Taking $k_1=0$ and $k_2=4$, by (\ref{y6}) and (\ref{y7}) we obtain
 \beq
 r_{ij}&=&\frac{1}{1+3b^2}(b_is_j+b_js_i),\label{y9}\\
F(\alpha,\beta)&=&(\alpha^2+4\beta^2)^{\frac{1}{4}}\sqrt{2\beta+\sqrt{\alpha^2+4\beta^2}}.\label{y11}
 \eeq
Theorem \ref{th1}(iv) shows that the metric (\ref{y11}) in
two-dimensional case is of isotropic S-curvature if and only if
$\beta$ satisfies (\ref{y9}). In the following example, we show a
pair $\alpha$ and $\beta$ such that (\ref{y9}) holds. For more
examples, see a general case in Example \ref{ex51} below in the
last section.

\begin{ex}\label{ex11}
 Let $F$ be a two-dimensional $(\alpha,\beta)$-metric defined by
 (\ref{y11}). Define $\alpha$ and $\beta$ by
 $$
 \alpha=e^{\sigma}\sqrt{(y^1)^2+(y^2)^2}, \ \ \ \beta=e^{\sigma}(\xi y^1+\eta
 y^2),
 $$
 where $\xi,\eta$ and $\sigma$ are scalar functions which are
 given by
  $$\xi=x^2, \ \ \eta=-x^1, \ \
  \sigma=-\frac{1}{4}\ln(1+4|x|^2), \ \ \ (|x|^2:=(x^1)^2+(x^2)^2).
  $$
  Then $\alpha$ and $\beta$ satisfy (\ref{y9}), and therefore, $F$
  is of isotropic S-curvature ${\bf S}=0$ by Theorem
  \ref{th1}(iv). Further we have $b^2=||\beta||_{\alpha}^2=|x|^2\ne constant$.
\end{ex}

\begin{rem}
 By a result in \cite{CMS}, the Finsler metric $F$ in Theorem
 \ref{th1} (iv) is an Einstein metric (equivalently, isotropic
 flag curvature), since ${\bf S}=0$.  We show further in
 \cite{Y3} that, the local structure of $F$ in Theorem
 \ref{th1} (iv) can be almost determined, and $F$ is an Einstein
 metric but generally not Ricci-flat (non-zero isotropic
 flag curvature). If taking
$k_1=-1,\ k_2=0,$
 then $F$ in Theorem
 \ref{th1} (iv) becomes $F=\sqrt{\alpha(\alpha+\beta)}$, which is
 called a square-root metric (\cite{Y2}). We show in \cite{Y2} that, a
 2-dimensional square-root metric $F$ is an Einstein metric if and
 only if $F$ is of vanishing S-curvature.
\end{rem}

\section{Preliminaries}

 Let $F$ be a Finsler metric on an $n$-dimensional manifold $M$
 with $(x^i,y^i)$  the standard local coordinate system in
 $TM$. The Hausdorff-Busemann volume form
$dV=\sigma_F(x)dx^1\wedge ... \wedge dx^n$ is defined by
  $$\sigma_F(x):=\frac{Vol(B^n)}{Vol\big\{(y^i)\in
  R^n|F(y^i\pxi|_x)<1\big \}}.$$
The Finsler metric $F$ induces a vector field $G=y^i\pxi-2G^i\pyi$
on $TM$ defined by
$$G^i=\frac{1}{4}g^{il}\big\{[F^2]_{x^ky^l}y^k-[F^2]_{x^l}\big\}.$$
Then the S-curvature is defined by
$${\bf S}:=\frac{\pa G^m}{\pa y^m}-y^m\frac{\pa}{\pa
x^m}(ln\sigma_F).$$
 ${\bf S}$ is said to be {\it isotropic} if there is a scalar function
 $c(x)$ on $M$ such that
 $${\bf S}=(n+1)c(x)F.$$
 If $c(x)$ is a constant, then we call $F$ is of {\it constant
 S-curvature}.

An $(\alpha,\beta)$-metric is expressed in the following form:
 $$F=\alpha \phi(s),\ \ s=\beta/\alpha,$$
where $\phi(s)>0$ is a $C^{\infty}$ function on $(-b_o,b_o)$. It
is known that $F$ is a positive definite regular Finsler metric
with $ \|\beta\|_{\alpha} < b_o$ if
 $$
 \phi(s)-s\phi'(s)+(\rho^2-s^2)\phi''(s)>0, \ \ (|s| \leq \rho
 <b_o).
 $$

For an $n$-dimensional $(\alpha,\beta)$-metric $F=\alpha \phi(s),
s=\beta/\alpha$, it has been shown in \cite{CS1} that the
S-curvature is given by
 \be\label{y13}
{\bf
S}=\big\{2\Psi-\frac{f'(b)}{bf(b)}\big\}(r_0+s_0)-\alpha^{-1}\frac{\Phi}{2\Delta^2}(r_{00}-2\alpha
Qs_0),
 \ee
where $\Phi$ is defined by (\ref{y01}) and
$$r_0:=r_iy^i,\ \ s_0:=s_iy^i, \ \ r_{00}:=r_{ij}y^iy^j,$$
$$
 \Psi:=\frac{Q'}{2\Delta},\ \ \Delta:=1+sQ+(b^2-s^2)Q', \ \
 Q:=\frac{\phi'}{\phi-s\phi'},
$$
\be\label{y14}
 f(b):=\frac{\int^{\pi}_0 sin^{n-2}tdt}{\int^{\pi}_0\frac{sin^{n-2}t}{\phi(b cost)^n}dt}.
\ee

Fix an arbitrary point $x\in M$ and take  an orthonormal basis
  $\{e_i\}$ at $x$ such that
   $$\alpha=\sqrt{\sum_{i=1}^n(y^i)^2},\ \ \beta=by^1.$$
Then we change coordinates $(y^i)$ to $(s, y^A)$ such that
  $$\alpha=\frac{b}{\sqrt{b^2-s^2}}\bar{\alpha},\ \
  \beta=\frac{bs}{\sqrt{b^2-s^2}}\bar{\alpha}, $$
where $\bar{\alpha}=\sqrt{\sum_{A=2}^n(y^A)^2}$. Let
 $$\bar{r}_{10}:=\sum^n_{A=2}r_{1A}y^A, \ \ \bar{r}_{00}:=\sum^n_{A,B=2}r_{AB}y^Ay^B, \ \
 \bar{s}_0:=\sum^n_{A=2}s_Ay^A.$$
 By (\ref{y13}), it is shown in \cite{CS1} that ${\bf
 S}=(n+1)c(x)F$ is equivalent to the following two equations:
  \be\label{y15}
  \frac{\Phi}{2\Delta^2}(b^2-s^2)\bar{r}_{00}=-\Big\{s\big[\frac{s\Phi}{2\Delta^2}-2\Psi
  b^2+\frac{bf'(b)}{f(b)}\big]r_{11}+(n+1)cb^2\phi\Big\}\bar{\alpha}^2,
  \ee
 \be\label{y16}
\Big\{\frac{s\Phi}{\Delta^2}-2\Psi
  b^2+\frac{bf'(b)}{f(b)}\Big\}r_{1A}=\Big\{(\frac{\Phi
  Q}{\Delta^2}+2\Psi)b^2-\frac{bf'(b)}{f(b)}\Big\}s_{1A}.
 \ee

It studies (\ref{y15}) and (\ref{y16}) in \cite{CS1} by three
steps: (i) $\Phi=0$, (ii) $\Phi \ne 0$ and $\Upsilon =0$, and
(iii) $\Phi \ne 0$ and $\Upsilon\ne 0$, where $\Upsilon$ is
defined by
 $$\Upsilon:=\frac{d}{ds}\Big[\frac{s\Phi}{\Delta^2}-2\Psi
 b^2\Big].$$

In the discussion for the third case $\Phi \ne 0$ and $\Upsilon\ne
0$, it obtains in \cite{CS1} the following lemma (see Lemma 6.1 in
\cite{CS1}):

\begin{lem} (\cite{CS1})
 Let $F=\alpha\phi(s),s=\beta/\alpha$, be an
 $(\alpha,\beta)$-metric on an $n$-dimensional manifold. Assume
 $\phi(s)$ satisfies $\Phi \ne 0$ and $\Upsilon\ne 0$, and $F$ has
 isotropic S-curvature, ${\bf
 S}=(n+1)c(x)F$.  Then
 \be\label{y17}
 r_{ij}=ka_{ij}-\epsilon b_ib_j-\lambda (b_is_j+b_js_i),
 \ee
 \be\label{y18}
 -2s(k-\epsilon b^2)\Psi+(k-\epsilon
 s^2)\frac{\Phi}{2\Delta^2}+(n+1)c\phi-s \nu=0,
 \ee
 where $\lambda=\lambda(x),k=k(x)$ and $\epsilon=\epsilon(x)$ are
 some scalar functions and
  \be\label{y19}
 \nu:=-\frac{f'(b)}{bf(b)}(k-\epsilon b^2).
  \ee
  If in addition $s_0\ne 0$, then
   \be\label{y20}
  -2\Psi-\frac{Q\Phi}{\Delta^2}-\lambda
  \big(\frac{s\Phi}{\Delta^2}-2\Psi b^2\big)=\delta,
   \ee
   where
   \be\label{y21}
\delta:=-\frac{f'(b)}{bf(b)}(1-\lambda b^2).
   \ee
\end{lem}

For the proof of Lemma 6.2 in \cite{CS1}, it uses  the method of
expressing (\ref{y18}) and (\ref{y20}) as polynomials of $b$ when
$b\ne constant$, but no consideration is taken on the effect of
$\lambda,k,c$, etc. We should point out that these functions are
actually dependent on $b$. Therefore, the method we adopt is to
expand (\ref{y18}) and (\ref{y20}) as power series of $s$.

\section{On equation (\ref{y18})}\label{sec3}
In this section, we assume $b\ne constant$ and  $\phi(s)\ne
k_1\sqrt{1+k_2s^2}+k_3s$ for any constants $k_1>0, k_2$ and $k_3$.
We are going to  prove that $k=0,c=0,\epsilon=0$ and $\nu=0$ in
(\ref{y18}).

We first transform (\ref{y18}) into a differential equation about
$\phi(s)$ and then (\ref{y18})$\times 2\phi
[\phi-s\phi'+(b^2-s^2)\phi'']^2$ yields
 \be
\Gamma_0=0,\label{y22}
 \ee
 where we omit the expression of $\Gamma_0$.

Let $p_i$ be the coefficients of $s^i$ in (\ref{y22}). We need to
compute $p_0,p_1,p_2,p_3$ and $p_4$ first. For this, it is
sufficient to plug
$$\phi(s)=1+a_1s+a_2s^2+a_3s^3+a_4s^4+a_5s^5+a_6s^6+a_7s^7+o(s^7)$$
into (\ref{y22}). Here we omit the expressions of
$p_0,p_1,p_2,p_3$ and $p_4$. All the equations $p_i=0$ are linear
equations about $k,c,\epsilon$ and $\nu$.  By
$p_0=0,p_1=0,p_2=0,p_3=0$ and $b\ne constant$, it is easy to
conclude that $k=0,c=0,\epsilon=0$ and $\nu=0$ if there hold
$$a_1\ne 0,\ \ a_4\ne -\frac{2(n+1)a_2^2+(n-2)a_1a_3}{4(n+1)},$$
since in this case the coefficient determinant of the linear
system $p_0=0,p_1=0,p_2=0,p_3=0$ is not zero.

In the following we prove there also hold $k=0,c=0,\epsilon=0$ and
$\nu=0$ if $a_1=0$, or $4(n+1)a_4+2(n+1)a_2^2+(n-2)a_1a_3=0$.

\bigskip

\noindent{\bf Case 1:} Assume $a_1=0$.  By $p_0=0,p_1=0$ and
$a_1=0$, we obtain
 \be\label{y24}
\nu=\frac{2\big[(18a_3^2-10a_2^3-12a_2a_4)b^4-(7a_2^2+6a_4)b^2-a_2\big]k+2a_2b^2(1+2a_2b^2)^2\epsilon}{(1+2a_2b^2)^3},
 \ee
 \be\label{y25}
 c=\frac{3a_3b^2}{(n+1)(1+2a_2b^2)^2}k.
 \ee
Since $\phi(s)\ne \sqrt{1+2a_2s^2}$, there exists some minimal
integer $m$ such that
 \be\label{y23}
  a_{2m+1}\ne 0,\ \ (m\ge 1); \ \ {\text or} \ \
 a_{2m}\ne C^m_{\frac{1}{2}}(2a_2)^m,\ \ (m\ge 2),
 \ee
 where $C^i_{\mu}$ are the
 generalized combination coefficients.

 \

{\bf Case 1A.} Assume $a_{2m+1}\ne 0$ in (\ref{y23}). If $a_3\ne
0$ ($m=1$), plug (\ref{y24}), (\ref{y25}) and $a_1=0$ into $p_2=0$
and  $p_4=0$ and then we get a linear system about $k$ and
$\epsilon$. The critical component of the determinant for this
linear system is given by
 $$(...)b^8+(...)b^6+(...)b^4+(...)b^2-3(n-1)(n+3)a_3^2,$$
where the omitted terms in the brackets of the above are all
constants. Now it is easy to see that if $a_3\ne 0$, then $k=0$
and $\epsilon=0$.
 Thus by (\ref{y24}) and (\ref{y25}) we
 have $c=0$ and $\nu=0$.

If $m>1$, plug (\ref{y24}), (\ref{y25}), $a_1=0$ and  $a_3=0$ into
(\ref{y22})$\times (1+2a_2b^2)^2$ and we get another new
differential equation about $\phi(s)$, which is denoted by
 \be
  \Gamma_{00}=0,\label{y26}
 \ee
where we omit the expressions of $\Gamma_{00}$ for simplicity.
Let $q_i$ be the coefficients of $s^i$ in
 (\ref{y26}). For our purpose to prove $k=0$ and $\epsilon=0$, we only need to
 compute $q_{2m-2}$ and $q_{2m}$. Express $\phi(s)$ as
  \be\label{y27}
 \phi(s)=g(s)+h(s),
  \ee
where
$$g(s):=1+\sum_{i=1}^{\infty}a_{2i}s^{2i},\ \
h(s):=\sum_{i=m}^{\infty}a_{2i+1}s^{2i+1}.$$
 By a simple analysis on (\ref{y26}), it is easy to compute $q_{2m-2}$ and
 $q_{2m}$. To get $q_{2m-2}$, it is sufficient to put
 $$g(s)=1+o(s),\ \ h(s)=a_{2m+1}s^{2m+1}+o(s^{2m+2}),$$
 and plug (\ref{y27}) into (\ref{y26}). Then by $q_{2m-2}=0$ we obtain
  \be\label{y28}
 -2m(4m^2-1)b^2a_{2m+1}(1+2a_2b^2)^2k=0.
  \ee
By (\ref{y28}) we have $k=0$. To get $q_{2m}$, it is sufficient to
put
 $$g(s)=1+a_2s^2+o(s^3),\ \ h(s)=a_{2m+1}s^{2m+1}+a_{2m+3}s^{2m+3}+o(s^{2m+4}),$$
 and plug (\ref{y27}) into (\ref{y26}). Then by $k=0$ and $q_{2m}=0$  we obtain
 \be\label{y29}
 2m(2m+1)^2a_{2m+1}b^2(1+2a_2b^2)^2\epsilon=0.
 \ee
By (\ref{y29}) we have $\epsilon=0$. Thus by (\ref{y24}) and
(\ref{y25}) we
 have $c=0$ and $\nu=0$.

 \

{\bf Case 1B.} Assume all $a_{2i+1}=0$ ($i\ge 0$), and assume
$a_{2m}\ne C^m_{\frac{1}{2}}(2a_2)^m$ in (\ref{y23}). In this
case, we may express $\phi(s)$ as
 \be\label{y30}
 \phi(s)=\sqrt{1+2a_2s^2}+\eta(s),
 \ee
where
 $$\eta(s):=\sum_{i=m}^{\infty}d_{2i}s^{2i},\ \ d_{2m}\ne 0.$$

If $m=2$,  plug (\ref{y24}), (\ref{y25}), $a_1=0$  and $a_3=0$
into $p_3=0$ and $p_5=0$ and then we get a linear system about $k$
and $\epsilon$. The critical component of the determinant for this
linear system is given by
 $$(...)b^4+(...)b^2-(n+1)(n+4)(2a_4+a_2^2)^2,$$
where the omitted terms in the brackets of the above are all
constants. Now it is easy to see that if $2a_4+a_2^2\ne 0$, then
$k=0$ and $\epsilon=0$.
 Thus by (\ref{y24}) and (\ref{y25}) we
 have $c=0$ and $\nu=0$.

If $m>2$, plug (\ref{y30}) and $a_4=-a_2^2/2$ into
(\ref{y26})$\times (1+2a_2s^2)^3$ and we obtain a differential
equation about $\eta(s)$, which is denoted by
 \be\label{y33}
 \Gamma_1=0.
 \ee
Here we omit the expression of $\Gamma_1$. Let $q_i$ be the
coefficients of $s^i$ in
 (\ref{y33}). For our purpose to prove $k=0$ and $\epsilon=0$, we only need to
 compute $q_{2m-3}$ and $q_{2m-1}$.

To get $q_{2m-3}$, it is sufficient to plug
$$\eta(s)=d_{2m}s^{2m}+o(s^{2m+1})$$
into (\ref{y33}). Then by $q_{2m-3}=0$ we have
 \be\label{y34}
 16m(4m^2-1)b^2(1+2a_2b^2)^2d_{2m}k=0.
 \ee
 By (\ref{y34}) we get $k=0$. To get $q_{2m-1}$, it is sufficient to plug
$$\eta(s)=d_{2m}s^{2m}+d_{2m+2}s^{2m+2}+o(s^{2m+3})$$
into (\ref{y33}). Then by $k=0$ and  $q_{2m-1}=0$ we have
 \be\label{y35}
 16m^2(2m-1)b^2(1+2a_2b^2)^2d_{2m}\epsilon=0.
 \ee
By (\ref{y35}) we get $\epsilon=0$. Thus by (\ref{y24}) and
(\ref{y25}) we
 have $c=0$ and $\nu=0$.

\bigskip

\noindent{\bf Case 2:} Assume
$$a_1\ne 0,\ \
4(n+1)a_4+2(n+1)a_2^2+(n-2)a_1a_3=0.$$

By a simple analysis on the coefficient determinant of the linear
system $p_0=0,p_1=0,p_2=0,p_3=0$, it is enough for us to prove
$k=0,c=0,\epsilon=0$ and $\nu=0$ under one of the following two
conditions
 \be\label{y36}
 a_3=0,\ \ a_4=-\frac{1}{2}a_2^2,\ \
 a_6=\frac{1}{6}[(n-2)a_1a_5+3a_2^3],
 \ee
and
 \be\label{y37}
  a_3=-\frac{(4n^3+15n^2+16)a_1^3}{36(n^2-1)},\ \ a_4= \frac{2(n+1)a_2^2+(n-2)a_1a_3}{4(n+1)},
 \ee
 \be\label{y38}
 a_5=\frac{(n+4)(4n^2-n+4)}{1440(n+1)^3(1-n)}T_0,\ \ a_6=\frac{T}{60(n+1)^2},
 \ee
where
 \beqn
 T_0:&=&a_1^3\big[2a_1^2n^3+5(3a_1^2-16a_2)n^2+(6a_1^2-160a_2)n+20(a_1^2-4a_2)\big],\\
 T:&=&a_1(10a_5+20a_2a_3-3a_1^2a_3)n^3+(30a_2^3-120a_3^2+45a_1a_2a_3-6a_1^3a_3)n^2+\\
 &&(60a_2^3+15a_1^3a_3-30a_1a_5-276a_3^2-105a_1a_2a_3)n+18a_1^3a_3-130a_1a_2a_3\\
 &&-48a_3^2+30a_2^3-20a_1a_5.
 \eeqn

{\bf Case 2A.} Assume (\ref{y36}).  Solving $p_0=0,p_1=0,p_2=0$
and $p_4=0$
 yields
 \be\label{y40}
 k=\frac{2(1+2a_2b^2)c}{a_1},\ \
 \epsilon=\frac{2(a_1^2-2a_2)(1+2a_2b^2)c}{a_1},
 \ee
 \be\label{y41}
   a_5=0,\ \
   \nu=\frac{2[(1+n+2a_2b^2)a_1^2-2a_2(1+2a_2b^2)]c}{a_1}.
 \ee

Plug (\ref{y40}) and (\ref{y41}) into (\ref{y22}) and then we get
 \be\label{y42}
 c( f_0+f_2b^2+f_4b^4)=0,
  \ee
where $f_0,f_2,f_4$ are some ODEs about $\phi(s)$, where we omit
the expressions of $f_0,f_2,f_4$. If $c\ne 0$, then by
(\ref{y42}), solving $f_0=0,f_2=0,f_4=0$ with $\phi(0)=1$ yields
$\phi(s)=a_1s+\sqrt{1+2a_2s^2}$. This case is excluded. So $c=0$.
Then by (\ref{y40}) and (\ref{y41}) we get $k=0,\epsilon=0,\nu=0$.

\

{\bf Case 2B.} Assume (\ref{y37})-(\ref{y38}). Plug (\ref{y37})
and (\ref{y38}) into $p_0=0,p_1=0,p_2=0$ and $p_4=0$ and we obtain
$k=0,\epsilon=0,\nu=0,c=0$.

\section{On equation (\ref{y20})}\label{sec4}
In this section, we assume $b\ne constant$ and  $\phi(s)\ne
k_1\sqrt{1+k_2s^2}+k_3s$ for any constants $k_1>0, k_2$ and $k_3$.
We are going to  show that (\ref{y20}) has non-trivial solutions
in case of dimension $n=2$.

We first transform (\ref{y20}) into a differential equation about
$\phi(s)$ and then (\ref{y20})$\times \phi (-\phi+s\phi')
[\phi-s\phi'+(b^2-s^2)\phi'']^2$ gives
 \be
 \Gamma_2=0,\label{y43}
 \ee
where we omit the expression of $\Gamma_2$.

 Let $p_i$ be the coefficients of $s^i$ in (\ref{y43}). We need
to compute $p_0,p_1,p_2$ and $p_3$ first, where we do not write
out their expressions. For this, it is sufficient to plug
 $$\phi(s)=1+a_1s+a_2s^2+a_3s^3+a_4s^4+a_5s^5+a_6s^6+o(s^6)$$
into (\ref{y43}).  In the following we will solve $\lambda$ and
$\delta$ in two cases.

\bigskip

\noindent{\bf Case 1:} Assume $a_1=0$ and $a_3=0$. We are going to
show that this case is excluded.

Plugging  $a_1=0$ and $a_3=0$ into  $p_0=0$ yields
 \be\label{y44}
 \delta=\frac{2a_2}{1+2a_2b^2}(\lambda b^2-1).
 \ee
Since $\phi(s)\ne \sqrt{1+2a_2s^2}$, there exists some minimal
integer $m$
  such that
 \be\label{y45}
  a_{2m+1}\ne 0,\ \ (m\ge 2); \ \ {\text or} \ \
 a_{2m}\ne C^m_{\frac{1}{2}}(2a_2)^m,\ \ (m\ge 2),
 \ee
 where $C^i_{\mu}$ are the
 generalized combination coefficients. Then we will determine $\lambda$ in the  two cases of (\ref{y45}).

\

{\bf Case 1A.} Assume $a_{2m+1}\ne 0$ in (\ref{y45}). In this
case, we need to compute $p_{2m-1}$. For this, express $\phi(s)$
as
  \be\label{y46}
 \phi(s)=g(s)+h(s),
  \ee
where
$$g(s):=1+\sum_{i=1}^{\infty}a_{2i}s^{2i},\ \
h(s):=\sum_{i=m}^{\infty}a_{2i+1}s^{2i+1}.$$
 By a simple analysis
on (\ref{y43}), to get $p_{2m-1}$, it is sufficient to put
 $$g(s)=1+2a_2s^2+o(s^3),\ \ h(s)=a_{2m+1}s^{2m+1}+o(s^{2m+2}),$$
 and plug (\ref{y46}) into (\ref{y43}). Then by $p_{2m-1}=0$, $a_{2m+1}\ne 0$ and (\ref{y44}) we obtain

 \be\label{y48}
 \lambda=\frac{1-2(2m-1)a_2b^2}{2mb^2}.
 \ee

{\bf Case 1B.} Assume all $a_{2i+1}=0$ ($i\ge 0$), and assume
$a_{2m}\ne C^m_{\frac{1}{2}}(2a_2)^m$ in (\ref{y45}).  Express
$\phi(s)$ as
 \be\label{y49}
 \phi(s)=\sqrt{1+2a_2s^2}+\eta(s),
 \ee
where
 $$\eta(s):=\sum_{i=m}^{\infty}d_{2i}s^{2i},\ \ d_{2m}\ne 0.$$

Plug (\ref{y49}) and (\ref{y44}) into
(\ref{y43})$\times(1+2a_2s^2)^{7/2}$ and we obtain a differential
equation about $\eta(s)$, which is denoted by
 \be\label{y50}
 \Gamma_3=0.
 \ee
Here we omit the expression of $\Gamma_3$. Let $q_i$ be the
coefficients of $s^i$ in (\ref{y50}). We need to compute
$q_{2m-2}$. For this, it is sufficient to put
 $$\eta(s)=d_{2m}s^{2m}+o(s^{2m+1})$$
and plug (\ref{y49}) into  (\ref{y50}). Then by $q_{2m-2}=0$, we
obtain
 \be\label{y51}
 \lambda=\frac{1-4(m-1)a_2b^2}{(2m-1)b^2}.
 \ee

\vspace{0.3cm}

Now we have solved $\lambda$ in the two cases of (\ref{y45}). It
is easy to see that (\ref{y48}) and (\ref{y51}) can be written in
the following form
 \be\label{y52}
 \lambda=\frac{1-2(k-1)a_2b^2}{kb^2},
 \ee
 where $k\ge 3$ is an integer.

 Plugging (\ref{y44}) and (\ref{y52}) into (\ref{y43}) yields
 \be\label{y53}
 f_0+f_2b^2+f_4b^4=0,
 \ee
where $f_0,f_2,f_4$ are some ODEs about $\phi(s)$, where we omit
the expressions of $f_0,f_2,f_4$. Then by (\ref{y53}), solving
$f_0=0,f_2=0,f_4=0$ with $\phi(0)=1$ yields
$$\phi(s)=\sqrt{1+2a_2s^2}.$$
This case is excluded.

\bigskip

\noindent{\bf Case 2:} Assume $a_1\ne 0$ or $a_3\ne 0$. We are
going to show that for one case, there are non-trivial solutions
for $\phi(s)$ in dimensions $n=2$.

\

{\bf Case 2A.} Assume $a_1=0$ and $a_3\ne 0$.  It follows
$a_4=-\frac{1}{2}a_2^2$ from $p_0=0,p_1=0,p_2=0$ and $a_1=0$. Then
by $p_0=0,p_1=0,p_3=0$, $a_1=0$ and $a_4=-\frac{1}{2}a_2^2$ we get
a contradiction.

\

{\bf Case 2B.} Assume $a_1\ne 0$. Solving $\lambda,\delta$ from
$p_0=0,p_1=0$ and then plugging them into $p_2=0$ yields
 \be\label{y55}
 a_4=-\frac{1}{2}a_2^2-a_1a_3,
 \ee
 \be\label{y56}
 a_5=-\frac{a_3\big[n^2a_1^3+(3a_3+20a_1a_2-6a_1^3)n+20a_1a_2-21a_3-7a_1^3\big]}{10(n+1)a_1},
 \ee
 \be\label{y57}
 (n-7)a_3^2(na_1^3+a_1^3-6a_3)=0.
 \ee

 By (\ref{y57}), we break our discussion into the following three steps.

 \

{\bf (I).} If $n=7$ and $a_3\ne 0$, solving $\lambda,\delta$ from
$p_0=0,p_1=0$ and then plugging them together with $n=7$,
(\ref{y55}) and (\ref{y56}) into $p_3=0$ yields a contradiction.

\

{\bf (II).} If $a_3=0$, then plug (\ref{y55}) and $a_3=0$ into
$p_0=0,p_1=0$ and we can get
 \be\label{y58}
 \lambda=a_1^2-2a_2,\ \
 \delta=\frac{na_1^2+(1+2a_2b^2)(a_1^2-2a_2)}{1+2a_2b^2}.
\ee
 Plugging (\ref{y58}) into (\ref{y43}) yields
 \be\label{y59}
 f_0+f_2b^2+f_4b^4=0,
 \ee
where $f_0,f_2,f_4$ are some ODEs about $\phi(s)$, where we omit
the expressions of $f_0,f_2,f_4$. Then by (\ref{y59}), solving
$f_0=0,f_2=0,f_4=0$ with $\phi(0)=1$ yields
$$\phi(s)=a_1s+\sqrt{1+2a_2s^2}.$$
This case is excluded.

\

{\bf (III).}  Assume
 \be\label{y60}
a_3=\frac{1}{6}(n+1)a_1^3.
 \ee
 Solving $\lambda,\delta$ from
$p_0=0,p_1=0$ and then plugging them together with
 (\ref{y55}), (\ref{y56}) and (\ref{y60}) into $p_3=0$
yields
 $$
 (...)b^2+(n+1)(n-2)a_1^4=0,
 $$
 which implies $n=2$. Solving $\lambda,\delta$ from
$p_0=0,p_1=0$ and then plugging them together with (\ref{y55}),
 (\ref{y60}) and $n=2$ into (\ref{y43}) yields
  \be\label{y61}
 f_0+f_2b^2+f_4b^4=0,
 \ee
where $f_0,f_2,f_4$ are some ODEs about $\phi(s)$ given by
 \beqn
      f_0:=\hspace{-0.5cm}&&\big[2(a_1^2-a_2)s(\phi-s\phi')+\phi'\big]s^2\phi\phi'''-s^2\big[1+(2a_2-3a_1^2)s^2\big]\phi(\phi'')^2+\\
  &&\Big\{(1-2a_2s^2)(\phi-s\phi')^2+\big[4+2(3a_1^2-4a_2)s^2\big]s\phi'(\phi-s\phi')+6s^2(\phi')^2\Big\}\phi''+\\
  &&\big[(3a_1^2-2a_2)(\phi-s\phi')^2+(4a_2-3a_1^2)s\phi'(\phi-s\phi')-3(\phi')^2\big](\phi-s\phi'),\\
      f_2:=\hspace{-0.5cm}&&\Big\{\big[(2a_2+a_1^2)(3a_1^2-2a_2)s^2+2(a_2-a_1^2)\big]s(\phi-s\phi')-(1-2a_2s^2)\phi'\Big\}\phi\phi'''\\
 &&\big[1-(2a_2+a_1^2)s^2\big]\big[1+(2a_2-3a_1^2)s^2\big]\phi(\phi'')^2+\Big\{\big[(2a_2+a_1^2)(3a_1^2-2a_2)s^2+\\
 &&4a_1^2\big](\phi-s\phi')^2+\big[4(2a_2+a_1^2)(3a_1^2-2a_2)s^2+2(6a_2-a_1^2)\big]s\phi'(\phi-s\phi')+\\
 &&3(4a_2s^2-1)(\phi')^2\Big\}\phi''+\Big\{(2a_2+a_1^2)(2a_2-3a_1^2)(3s\phi'-\phi)(\phi-s\phi')\\
 &&-6a_2(\phi')^2\Big\}(\phi-s\phi'),\\
     f_4:=\hspace{-0.5cm}&&\big[(2a_2+a_1^2)(2a_2-3a_1^2)s(\phi-s\phi')-2a_2\phi'\big]\phi\phi'''+\\
 &&(2a_2+a_1^2)\big[1+(2a_2-3a_1^2)s^2\big]\phi(\phi'')^2+\\
 &&
 \big[(2a_2+a_1^2)(2a_2-3a_1^2)(\phi-s\phi')(3s\phi'-\phi)-6a_2(\phi')^2\big]\phi''.
 \eeqn
Then by (\ref{y61}), solving $f_0=0,f_2=0,f_4=0$ with $\phi(0)=1$
yields (\ref{y7}), where we put
 \be\label{k12}
 k_1:=2a_2-3a_1^2,\ \ k_2:=2a_2+a_1^2.
 \ee

Plugging (\ref{y55}) and
 (\ref{y60}) and $n=2$ into $p_0=0,p_1=0$ yields
  \be\label{y62}
 \delta=\frac{(3a_1^2-2a_2)\big[1+(2a_2+a_1^2)b^2\big]}{1+2a_2b^2},
  \ee
  \be\label{y63}
 \lambda=\frac{(3a_1^2-2a_2)(2a_2+a_1^2)b^2+2(a_1^2-a_2)}{1+2a_2b^2}.
  \ee
Since we have proved in Section 3 that $k=0,\epsilon=0$, by
(\ref{y17}) and (\ref{y63}) we obtain (\ref{y6}). By (\ref{y21}),
(\ref{y62}) and (\ref{y63}) we get
 \be\label{y64}
 f(b)=\sqrt{1+(2a_2-3a_1^2)b^2}.
 \ee

One possibly wonders whether we can get (\ref{y64}) from
(\ref{y14}) when we plug (\ref{y7}) and $n=2$ into (\ref{y14})?
This is true. One way to check it is to expand (\ref{y14}) and
(\ref{y64}) into power series respectively. One may try a direct
verification. We do not need this fact in our proof of Theorem
\ref{th1}.

\section{Further characterizations  and examples}
In this section,  we will discuss, for dimension $n=2$, the
solutions to the equations (\ref{y1}), (\ref{y3}), (\ref{y5}) and
(\ref{y6}) respectively, and construct some corresponding examples
for every class in Theorem \ref{th1}.

Since every two-dimensional Riemann metric is locally conformally
 flat, we may put
  \be\label{y69}
  \alpha=e^{\sigma}\sqrt{(y^1)^2+(y^2)^2},
  \ee
 where $\sigma=\sigma(x)$ is a scalar function and $x=(x^1,x^2)$. Then $b=constant$ is equivalent to
  \be\label{y70}
 \beta=\frac{be^{\sigma}(\xi y^1+\eta y^2)}{\sqrt{\xi^2+\eta^2}},
  \ee
where $\eta=\eta(x)$ and $\xi=\xi(x)$ are  scalar functions. If
$b\ne constant$, then $\beta$ can be expressed as
 \be\label{y71}
 \beta=e^{\sigma}(\xi y^1+\eta y^2).
 \ee

Note that actually (\ref{y69}) and (\ref{y70}) give all the local
solutions of (\ref{y1}).

\begin{prop}
 Let $F=\alpha \phi(s)$, $s=\beta/\alpha$, be a two-dimensional
  $(\alpha,\beta)$-metric on $R^2$, where $\phi(s)\ne \sqrt{1+k_1s^2}+k_2s$ (any constant $k_1,k_2$) satisfies
  (\ref{y2}) with $\phi(0)=1$. Then $F$ is of isotropic S-curvature  if and only if $\alpha$ and $\beta$
 can be locally defined  by (\ref{y69}) and (\ref{y70}), where
 $\xi,\eta$ and $\sigma$ are some scalar functions. In this case,
 ${\bf S}=0$.
\end{prop}

We consider
 \be\label{y72}
  r_{ij}=\epsilon(b^2a_{ij}-b_ib_j)-\frac{1}{b^2}(b_is_j+b_js_i),
 \ee
where $\epsilon=\epsilon(x)$ is a scalar function. It it easy to
show that in two-dimensional case, $b=constant$ if and only if
(\ref{y72}) holds. Assume  $\alpha$ and $\beta$
 are locally defined by (\ref{y69}) and (\ref{y70}) satisfying (\ref{y72}), then we have
\be\label{y73}
 \epsilon=\frac{(\xi^2+\eta^2)(\xi\sigma_1+\eta\sigma_2)-\xi\eta\eta_1+\xi^2\eta_2
 +\eta^2\xi_1-\xi\eta\xi_2}{be^{\sigma}(\xi^2+\eta^2)^{\frac{3}{2}}},
 \ee
where $\sigma_1:=\pa \sigma/\pa x^1$, $\sigma_2:=\pa \sigma/\pa
x^2$, etc.  Further, $\beta$ is  closed if and only if
 \be\label{y74}
(\xi^2+\eta^2)(\xi\sigma_2-\eta\sigma_1)-\xi^2\eta_1-\xi\eta\eta_2
 +\xi\eta\xi_1+\eta^2\xi_2= 0.
 \ee

Thus, it is easy to show that (\ref{y3}) can be locally
characterized by (\ref{y69}), (\ref{y70}) and (\ref{y74}), and
that (\ref{y5}) can be locally characterized by (\ref{y69}),
(\ref{y70}), (\ref{y74}), and $\epsilon=0$ in (\ref{y73}).

Suppose (\ref{y74}), and $\epsilon=0$ in (\ref{y73}), then we have
 \be\label{y75}
 \sigma_1=\frac{\eta\xi_2-\xi\eta_2}{\xi^2+\eta^2},\ \
 \sigma_2=\frac{\eta\xi_1-\xi\eta_1}{\xi^2+\eta^2}.
 \ee
We can find the  integrable condition of (\ref{y75}).

Put $\xi=x^2,\eta=1-x^1$ and plugging them into (\ref{y74}) yields
$$
 \sigma=f(\frac{x^2}{x^1-1})-ln(1-x^1),
 $$
where $f(t)$ is an arbitrary function. Further,  the above
$\sigma$ satisfies $\epsilon=0$ in (\ref{y73}) if and only if
$$
 f(t)=-\frac{1}{2}ln(1+t^2)+a,
$$
where $a$ is a constant.

\begin{prop}
Let $F=\alpha \phi(s)$, $s=\beta/\alpha$,  be a two-dimensional
  $(\alpha,\beta)$-metric on $R^2$, where $\phi(s)\ne \sqrt{1+k_1s^2}+k_2s$ (any constant $k_1,k_2$) satisfies
  (\ref{y4}) ($k\ne 0$) with $\phi(0)=1$. Then $F$ is of isotropic S-curvature ${\bf
  S}=3c(x)F$ if and only if $\alpha$ and $\beta$
 can be locally defined by (\ref{y69}) and (\ref{y70}), where
 $\xi,\eta$ and $\sigma$ are some scalar functions satisfying (\ref{y74}). In this case,
 ${\bf S}=3cF$ with $c=k\epsilon$, where $\epsilon$ is given
 by (\ref{y73}).
\end{prop}

\begin{prop}
Let $F=\alpha \phi(s)$, $s=\beta/\alpha$,  be a two-dimensional
  $(\alpha,\beta)$-metric on $R^2$, where $\phi(0)=1$, $b=constant$ and $\phi(s)\ne \sqrt{1+k_1s^2}+k_2s$ (any constant $k_1,k_2$) dose not
  satisfy (\ref{y4}) for any constant $k$. Then $F$ is of isotropic S-curvature if and only if $\alpha$ and $\beta$
 can be locally defined by (\ref{y69}) and (\ref{y70}), where
 $\xi,\eta$ and $\sigma$ are some scalar functions satisfying (\ref{y75}). In this case,
 ${\bf S}=0$.
\end{prop}

Finally we consider (\ref{y6}). In this case, $\alpha$ and $\beta$
can be expressed as (\ref{y69}) and (\ref{y71}). Let  $a_1,a_2$
satisfy (\ref{k12}) in the following.  We can show that (\ref{y6})
is equivalent to the following system of PDEs:
 \be\label{y76}
 \sigma_1=\frac{T_1}{T_0},\ \ \sigma_2=\frac{T_2}{\xi T_0},\ \
 \xi_1=-\frac{\eta(\eta\eta_2+\xi\xi_2+\xi\eta_1)}{\xi^2},
 \ee
where
 \beqn
 T_0:=\hspace{-0.5cm} &&\xi\big[1+(2a_2+a_1^2)(\xi^2+\eta^2)\big]\big[1+(2a_2-3a_1^2)(\xi^2+\eta^2)\big],\\
 T_1:=\hspace{-0.5cm}&&2\xi\eta\big[2(a_2-a_1^2)+(2a_2+a_1^2)(2a_2-3a_1^2)(\xi^2+\eta^2)\big]\xi_2\\
 && -\big[1+2(2a_2-a_1^2)\xi^2+2a_1^2\eta^2+(2a_2+a_1^2)(2a_2-3a_1^2)(\xi^4-\eta^4)\big]\eta_2,\\
 T_2:=\hspace{-0.5cm}&&\big[2a_1^2\xi^2+2(2a_2-a_1^2)\eta^2-(2a_2+a_1^2)(2a_2-3a_1^2)(\xi^4-\eta^4)\big](\xi\xi_2+\eta\eta_2)+\\
 &&
 \xi\big[1+(2a_2+a_1^2)(\xi^2+\eta^2)\big]\big[1+(2a_2-3a_1^2)(\xi^2+\eta^2)\big]\eta_1.
 \eeqn

\begin{prop}
Let $F=\alpha \phi(s)$, $s=\beta/\alpha$,  be a two-dimensional
  $(\alpha,\beta)$-metric on $R^2$, where  $b=||\beta||_{\alpha}\ne constant$ and $\phi(s)$
  satisfies (\ref{y7}). Then $F$ is of isotropic S-curvature  if and only if $\alpha$ and $\beta$
 can be locally defined by (\ref{y69}) and (\ref{y71}), where
 $\xi,\eta$ and $\sigma$ are some scalar functions satisfying (\ref{y76}). In this case,
 ${\bf S}=0$.
\end{prop}

If we take $\xi=x^2$ and $\eta=-x^1$, then $\sigma$ determined by
(\ref{y76}) is given by
 \be\label{y0067}
 \sigma=-\frac{1}{4}\Big\{ln\big[1+(2a_2+a_1^2)|x|^2\big]+3ln\big[1+(2a_2-3a_1^2)|x|^2\big]\Big\},
 \ee
where $|x|^2:=(x^1)^2+(x^2)^2$. Thus we obtain the following
example satisfying Theorem \ref{th1} (iv).

\begin{ex}\label{ex51}
 Let $F$ be a two-dimensional $(\alpha,\beta)$-metric defined by
 (\ref{y7}). Define $\alpha$ and $\beta$ by (\ref{y69}) and
 (\ref{y71}),
 where  $\xi=x^2$ and $\eta=-x^1$,
  and $\sigma$ is given by (\ref{y0067}).
  Then  $F$
  is of isotropic S-curvature ${\bf S}=0$ by Theorem
  \ref{th1}(iv). Further we have $b^2=||\beta||_{\alpha}^2=|x|^2\ne constant$.
\end{ex}

In Example \ref{ex51}, if we take $2a_2-3a_1^2=0$ and $a_1=1$,
then by (\ref{y7}) and (\ref{y8}), we obtain
 $$
\phi(s)=(1+4s^2)^{\frac{1}{4}}\sqrt{2s+\sqrt{1+4s^2}},
 $$
and thus we get Example \ref{ex11}.

\vspace{0.6cm}

\noindent Guojun Yang \\
Department of Mathematics \\
Sichuan University \\
Chengdu 610064, P. R. China \\
{\it  e-mail :} yangguojun@scu.edu.cn

\end{document}